# The combined equilibrium of business land-use and its congestion pricing principles (working paper)

Qian Liu [a,*], Chongchao Huang [a]

[a] School of Mathematics and Statistics, Wuhan University, Hubei, 430072, China

[*] Corresponding author.

E-mail addresses: mathlq@163.com (Q. Liu), cchuang@whu.edu.cn (C. Huang).

**Abstract**

This working paper is divided into two parts. Firstly, we develop a new combined equilibrium model of business land-use, which puts travelers' traffic equilibrium and business companies' competitive location equilibrium into a unified framework. A variational inequality is presented for the combined equilibrium and the properties of equilibrium solution are investigated. Secondly, the congestion pricing principles associated with the combined equilibrium are studied. From the mathematical point of view, we prove that there exists an optimal road pricing scheme that can minimize the social cost of travelers. This road pricing scheme generalizes the traditional link-based optimal road pricing scheme. Furthermore, when simultaneously imposing charges on travelers and companies is allowed, we prove that there exists an optimal congestion pricing scheme that can derive a combined equilibrium toward an overall system optimum. The economic meaning of every pricing scheme proposed in this paper is discussed in detail. At last, a simple numerical example is used to demonstrate that the optimal congestion pricing scheme may indeed reduce the social cost.

# 1. Introduction

Traffic congestion is one of the basic characteristics in the development of modern city. The development of modern economy derives a great amount of travel demand, which is often close to or even reaches the physical capacity of the whole traffic system and leads to broad traffic congestion. Recently, many metropolises are becoming more and more crowded and there is not enough space to expand the physical capacity of urban transportation system to reduce the congestion. Thus, drawing support from scientific management methods provides government with an alternative solution to mitigate traffic congestion. Congestion pricing, as a particular method, has been put into practice to effectively control traffic congestion. In academic field, there have been many studies on road congestion pricing. For a comprehensive review of classical mathematical theories and methods of road pricing, please see the monograph (Yang and Huang, 2005). The current researches about road pricing in the literature mainly focus on the more practical and complicated scheme design of road pricing (e.g. Hamdouch and Lawphongpanich, 2010; Yang, Xu, He and Meng, 2010; Lawphongpanich and Yin, 2012), the welfare effect of road pricing policy on social groups with different income levels (e.g. Wu, Yin and Lawphongpanich, 2011; Wu, Yin, Lawphongpanich and Yang, 2012), the congestion management with tradable credits (e.g. Yang and Wang, 2011; Wang et al., 2012; Wu et al., 2012) and so on.

A potential restriction of traditional road pricing policy to mitigate congestion is that it can not fundamentally alter the travel demand pattern in the network i.e. the travel demands between different O-D pairs in the network. The travel demand pattern is mainly determined by the land-use pattern while the road pricing policy can significantly change the travelers' travel route choices but only influence the number of travel demand slightly between a given O-D pair. Therefore, in order to fundamentally mitigate traffic congestion, we should explicitly consider the influence of the land-use pattern on the travel demand pattern of urban road network.

The close relationship between transport and land use has been recognized and studied for decades. For the research results in recently years, we refer to Briceño et al. (2008), Bravo et al. (2010), Meng et al. (2009), and Ma and Lo (2012), just to name a few. As far as we are concerned, almost all of the researches about the interactions between transport and land use in existing literature are limited to consider the behaviors of the group of travelers, though some classified the travelers into several kinds (e.g. Bravo et al., 2010). However, the trip distribution is not merely determined by travelers themselves but also clearly affected by the business location distribution i.e. the number of business companies located at every business center of the network. Indeed, if every destination of the network is supposed to be a business center, then the trip attraction at each destination to travelers (consumers) obviously increases with the number of business companies in this destination. From consumers' view, more business companies imply more shopping options and greater shopping convenience. Plenty of companies within a prosperous business center, for example, can always attract many travelers (consumers) even when the center itself and the paths leading to the center have been crowded. Additionally, it can be observed that a business center would attract many business companies to join with high rental rates if the number of travelers gathering in this center is large. The reason is that more consumers mean more potential business opportunities for business companies. Thus, the business location distribution pattern of companies and the trip distribution pattern of travelers can affect each other, and these effects should not be ignored. Furthermore, the business location distribution pattern is influenced by the congestion pricing policy as well. For example, if congestion tolls are set on the roads leading to a business center, then some previous customers of this business center may go to others. As a result, the original business competitive location equilibrium among different business centers in the network is broken, which will promote the business companies to relocate. This means that the change of business location distribution pattern, in the long term, will take place and influence the trip distribution of travelers in turn. Consequently, it makes sense to study the combined equilibrium of business land-use, which consists of business competitive location and transport system, and investigate corresponding congestion pricing principles. To the best of our knowledge, however, the mathematical analysis of the combined equilibrium of business land-use remains open, let alone the investigation of its congestion pricing principles.

The aim of this working paper is to investigate the combined equilibrium of business land-use by putting business competitive location equilibrium and traffic equilibrium into a unified framework, and derive corresponding optimal congestion pricing schemes. We first propose a new way of congestion pricing---business congestion tax policy, which imposes an extra tax on the companies in every business center of the network. Specifically, we charge the companies in the business centers which are too crowded or attract too many travelers (consumers) to result in network congestion, and provide subsidies for the companies in the business centers which attract relatively fewer travelers (consumers). This policy is to urge a portion of companies in original crowded business centers to relocate to other uncrowded business centers so as to balance the companies' distribution in different business centers. Thus every business center's attraction measure for travelers would change and fundamentally influence their destination choices. As a result, the traffic demands between different OD pairs would be balanced, the travelers would not gather at only a few business centers and the network congestion would be fundamentally mitigated.

The rest of the paper is organized as follows. In Section 2, we firstly give the definition of combined equilibrium of business land-use. The combined equilibrium consists of two sub-equilibriums---parametric traffic equilibrium and parametric business competitive location equilibrium which are dependent on each other. Two parametric variational inequalities are presented for these two sub-equilibrium problems. Integrating the two parametric variational inequalities together, we obtain a new variational inequality model equivalent to the combined equilibrium issue. With the help of this variational inequality, we confirm the existence and uniqueness of the combined equilibrium under some mild assumptions. In Section 3, we prove there is a road pricing scheme which can support a combined equilibrium as a traffic system optimum and give its mathematical formulas. Next, we describe the reasonability of the business congestion tax policy and prove that there exists an optimal congestion pricing scheme that includes both road pricing and business congestion tax to minimize the total social cost of travelers and business companies. The potential economic meaning of every pricing scheme is discussed. A simple numerical example is offered in Section 4 to demonstrate that the congestion pricing scheme consisting of road pricing and business congestion tax may indeed reduce the social cost. Finally, Section 5 summarizes and concludes the paper.

## 2. Model and mathematical formulations

2.1 Definition of combined equilibrium

The main goal of this section is to establish the mathematical models for the combined equilibrium of business land-use and prove its existence and uniqueness. It is assumed throughout the paper that every traveler in the network is a consumer who desires to choose the "optimal" destination and the fastest path leading to this "optimal" destination to go shopping. Now we list the terminologies and symbols used throughout this paper.

Consider a general network $G = (V, A)$, together with $N$ as node set and $A$ as directed link set. Let $O$, $D$ respectively denote the origin set and the destination set. For any $r \in O, s \in D$, let $R_{rs}$ denote the set of all simple routes between $OD$ pair $rs$. The total travel demand starting from origin node $r \in O$ is assumed to be fixed and given by $O_r$. Let $q = (q_{rs}, r \in O, s \in D)^T$ be the vector of travel demands where $q_{rs}$ is the travel demand between $OD$ pair $rs$, $d = (d_s, s \in D)^T$ be the vector of destination demands where $d_s = \sum_{r \in O} q_{rs}$ is the number of travelers choosing $s$ as their destination, $f = (f_k^{rs}, k \in R_{rs}, r \in O, s \in D)^T$ be the vector of path flows where $f_k^{rs}$ is the traffic flow on path $k \in R_{rs}$ and $x = (x_a, a \in A)^T$ be the vector of link flows where $x_a$ is the traffic link flow on directed link $a \in A$. The symbol $\delta_{ak}^{rs}$ equals 1 if link $a$ is on path $k \in R_{rs}$, otherwise equals 0. Let $t_a(x_a)$ be the separable travel time function of link $a$ which is dependence of link flow $x_a$. Throughout this paper, every destination node $s \in D$ is assumed to be a business center in the network. Let $h = (h_s, s \in D)^T$ denote the vector of business flows, where $h_s$ is the number of business companies locating in business center $s \in D$. We set that $T$ is the number of total business companies in the network, so $T = \sum_{s \in D} h_s$. In particular, $T$ is assumed to be a constant strictly greater than zero. For any destination node $s \in D$, we respectively denote $A_s(d_s, h_s)$ as the trip attraction function of $s$ for travelers and $B_s(d_s, h_s)$ as the business attraction function of $s$ for business companies. These two functions both depend on $d_s = \sum_{r \in O} q_{rs}$ and $h_s$.

Let $K_1$ be the set of all feasible trip distribution patterns defined by:

$$K_1 = \left\{ (f, x, q, d) \left| \begin{array}{l} f_k^{rs} \geq 0, \sum_{k \in R_{rs}} f_k^{rs} = q_{rs}, \sum_{s \in D} q_{rs} = O_r, \sum_{r \in O} \sum_{s \in D} \sum_{k \in R_{rs}} \delta_{ak}^{rs} f_k^{rs} = x_a, \sum_{r \in O} q_{rs} = d_s \\ k \in R_{rs}, r \in O, s \in D, a \in A \end{array} \right. \right\},$$

$K_2$ be the set of all feasible business location patterns defined by:

$$K_2 = \left\{ h \left| h_s \geq 0, \sum_{s \in D} h_s = T, s \in D \right. \right\}.$$

We introduce the following assumptions throughout this paper:

**Assumption 1.** For any link $a \in A$, the travel time function $t_a(x_a)$ is set to be separable, convex, continuous differentiable and strictly increasing of $x_a$. □

**Assumption 2.** For any $s \in D$, the destination attraction function for travelers, $A_s(d_s, h_s)$ is continuous differentiable of $(d_s, h_s)^T$, strictly decreasing of $d_s$, strictly increasing of $h_s$ and strictly concave with respect to $d_s$. □

**Assumption 3.** For any business center $s \in D$, the business attraction function for business companies, $B_s(d_s, h_s)$ is continuous differentiable of $(d_s, h_s)^T$, strictly increasing of $d_s$ and strictly decreasing of $h_s$. □

**Assumption 4.** Throughout this paper, we suppose that the functions $t_a(x_a)$, $A_s(d_s, h_s)$ and $B_s(d_s, h_s)$ have the same unit of measurement. That is to say, we do not distinguish between money-based cost and time-based cost in the paper. □

Now we formally introduce the definition of combined equilibrium.

**Definition 1.**

Let $(f, x, q, d) \in K_1$ and $h \in K_2$, we say $(f, x, q, d, h)$ is a (logit-based) combined equilibrium of business land-use if $(f, x, q, d, h)$ satisfies:

$$q_{rs} = O_r \frac{\exp\{-\gamma[v_{rs} - A_s(d_s, h_s)]\}}{\sum_{j \in D} \exp\{-\gamma[v_{rj} - A_j(d_j, h_j)]\}}, \quad \forall r \in O, s \in D \quad (1)$$

$$f_k^{rs} = q_{rs} \frac{\exp\{-\alpha[\sum_{a \in A} \delta_{ak}^{rs} t_a(x_a)]\}}{\sum_{j \in R_{rs}} \exp\{-\alpha[\sum_{a \in A} \delta_{aj}^{rs} t_a(x_a)]\}}, \quad \forall k \in R_{rs}, r \in O, s \in D; \quad (2)$$

$$h_s = T \frac{\exp[\beta B_s(d_s, h_s)]}{\sum_{j \in D} \exp[\beta B_j(d_j, h_j)]}, \quad \forall s \in D \quad (3)$$

where $\alpha, \beta, \gamma > 0$ and $v_{rs}$ is the logit-based expected travel time between OD pair $rs$, i.e. $v_{rs} = -\frac{1}{\alpha} \ln\{\sum_{j \in R_{rs}} \exp[-\alpha(\sum_{a \in A} \delta_{aj}^{rs} t_a(x_a))]\}$.

**Remark:** From the equations (1) and (2) in above definition, we know that the behaviors of travelers (consumers) have a nested logit choice structure. Additionally, by the logit choice theory, we have $\alpha \geq \gamma$.

In this paper, we define that $(f, x, q, d) \in K_1$ is the parametric traffic equilibrium if $(f, x, q, d)$ meets (1) and (2), where $h \in K_2$ is viewed as a parametric vector. Similarly, we define that $h \in K_2$ is the parametric business competitive location equilibrium if $h$ meets (3), where $(f, x, q, d) \in K_1$ is viewed as a parametric vector.

The goal of Section 2 is to establish the mathematical formulations for the combined equilibrium. Our strategy is to firstly divide the combined equilibrium into two interdependent parts: parametric traffic equilibrium and parametric business competitive location equilibrium and study them separately. Then we would integrate them together to derive the model of the combined equilibrium.

2.2 Parametric traffic equilibrium

Now we present a parametric variational inequality for the parametric traffic equilibrium problem.

**(VI 1)**

Find a vector $(f^*, x^*, q^*, d^*) \in K_1$ such that

$$\left[\begin{pmatrix} f \\ x \\ q \\ d \end{pmatrix} - \begin{pmatrix} f^* \\ x^* \\ q^* \\ d^* \end{pmatrix}\right]^T F(f^*, x^*, q^*, d^*, h) \geq 0, \quad \forall (f, x, q, d) \in K_1 \tag{4}$$

where $h = (h_s, s \in D)^T \in K_2$ is viewed as a parametric vector and

$$F(f, x, q, d, h) = \left[\frac{1}{\alpha}\ln f_k^{rs}, t_a(x_a), (\frac{1}{\gamma} - \frac{1}{\alpha})\ln q_{rs}, -A_s(d_s, h_s), k \in R_{rs}, r \in O, s \in D, a \in A\right]^T$$

.

The left term in (4) can be reformulated as:

$$\sum_{r \in O}\sum_{s \in D}\sum_{k \in R_{rs}}(f_k^{rs} - f_k^{rs*})\frac{1}{\alpha}\ln f_k^{rs*} + \sum_{a \in A}(x_a - x_a^*)t_a(x_a^*) + \sum_{r \in O}\sum_{s \in D}(q_{rs} - q_{rs}^*)(\frac{1}{\gamma} - \frac{1}{\alpha})\ln q_{rs}^*$$
$$-\sum_{s \in D}(d_s - d_s^*)A_s(d_s^*, h_s) \geq 0$$

We have the following result:

**Theorem 1.**

Let $h \in K_2$ and $(f,x,q,d) \in K_1$, then $(f,x,q,d)$ is a solution of (VI 1) if and only if $(f,x,q,d)$ is a parametric traffic equilibrium.

## 2.3 Parametric business competitive location equilibrium

Now we proceed to present another parametric variational inequality for the parametric business competitive location equilibrium.

**(VI 2)**

Find a vector $h^* = (h_s^*, s \in D)^T \in K_2$ such that

$$(h-h^*)^T U(d,h^*) \geq 0, \quad \forall h \in K_2, \tag{22}$$

where $U(d,h) = \left(\frac{1}{\beta} \ln h_s - B_s(d_s, h_s), s \in D\right)^T$, $(f,x,q,d)$ is a parametric vector.

The left term in (22) can be reformulated as: $\sum_{s \in D}(h_s - h_s^*)[\frac{1}{\beta} \ln h_s^* - B_s(d_s, h_s^*)] \geq 0$.

We have the following result.
**Theorem 2**.

Let $(f,x,q,d) \in K_1$ and $h \in K_2$, then $h$ is a solution of (VI 2) if and only if, it is a parametric business competitive location equilibrium.

## 2.4 Combined equilibrium of business land-use

So far we have studied the parametric traffic equilibrium and the parametric business competitive location equilibrium in a separate way. Now the aim is to put them into a framework to derive the model of combined equilibrium.

### 2.4.1 Variational inequalities formulation

In the following, a system of variational inequalities model (VIS 3), which is a combination of (VI 1) and (VI 2), is presented for the combined equilibrium.

**(VIS 3)**

Find a vector $(f^*, x^*, q^*, d^*, h^*) \in K_1 \times K_2$ such that

$$\left[\begin{pmatrix} f \\ x \\ q \\ d \end{pmatrix} - \begin{pmatrix} f^* \\ x^* \\ q^* \\ d^* \end{pmatrix}\right]^T F(f^*, x^*, q^*, d^*, h^*) \geq 0 \text{ and } (h-h^*)^T U(d^*, h^*) \geq 0,$$

$$\forall (f, x, q, d, h) \in K_1 \times K_2,$$

where $F(f^*, x^*, q^*, d^*, h^*)$ and $U(d^*, h^*)$ are respectively given in (4) and (22).

To discuss the solution existence of (VIS 3) conveniently, we need to reformulate (VIS 3) as a single variational inequality (VI 4) which is defined as follows:

**(VI 4)**

Find a vector $(f^*, x^*, q^*, d^*, h^*) \in K_1 \times K_2$ such that

$$\left[\begin{pmatrix} f \\ x \\ q \\ d \\ h \end{pmatrix} - \begin{pmatrix} f^* \\ x^* \\ q^* \\ d^* \\ h^* \end{pmatrix}\right]^T \psi(f^*, x^*, q^*, d^*, h^*) \geq 0, \quad \forall (f, x, q, d, h) \in K_1 \times K_2 \quad (23)$$

where $\psi(f, x, q, d, h) = \begin{pmatrix} F(f, x, q, d, h) \\ U(d, h) \end{pmatrix}$. The left term in (23) can be rewritten as

$$\sum_{r \in O} \sum_{s \in D} \sum_{k \in R_{rs}} (f_k^{rs} - f_k^{rs*}) \frac{1}{\alpha} \ln f_k^{rs*} + \sum_{a \in A} (x_a - x_a^*) t_a(x_a^*) + \sum_{r \in O} \sum_{s \in D} (q_{rs} - q_{rs}^*) \cdot (\frac{1}{\gamma} - \frac{1}{\alpha}) \ln q_{rs}^*$$

$$- \sum_{s \in D} (d_s - d_s^*) A_s(d_s^*, h_s^*) + \sum_{s \in D} (h_s - h_s^*) [\frac{1}{\beta} \ln h_s^* - B_s(d_s^*, h_s^*)] \geq 0$$

**Proposition 3.**

Suppose $(f^*, x^*, q^*, d^*) \in K_1$ and $h^* \in K_2$. Then $(f^*, x^*, q^*, d^*, h^*)$ is a solution of (VIS 3) if and only if it is a solution of (VI 4).

$\psi(f, x, q, d, h)$ is continuous on the compact convex set $K_1 \times K_2$, so (VI 4) has a solution. Therefore, the existence of the combined equilibrium can be confirmed.

2.4.2 Uniqueness of the combined equilibrium

Now we focus on the uniqueness of the combined equilibrium i.e. the solution uniqueness of variational inequality (VI 4). The uniqueness of the combined equilibrium pattern is significant to the study of optimal congestion pricing scheme. If the uniqueness can be confirmed, then the congestion pricing scheme's effect on trip distribution and business location distribution becomes predictable and some appropriate pricing scheme can be designed for achieving a target distribution pattern of travelers and business companies. To be precise, we have the following result.

**Theorem 4.**

The variational inequality (VI 4) has a unique solution if the following two conditions are valid for any $(f, x, q, d, h) \in K_1 \times K_2$:

(1) $\dfrac{\partial A_s(d_s, h_s)}{\partial d_s} + \dfrac{1}{2} \cdot \dfrac{\partial A_s(d_s, h_s)}{\partial h_s} + \dfrac{1}{2} \cdot \dfrac{\partial B_s(d_s, h_s)}{\partial d_s} < 0$, $\forall s \in D$;

(2) $\dfrac{\partial B_s(d_s, h_s)}{\partial h_s} + \dfrac{1}{2} \cdot \dfrac{\partial B_s(d_s, h_s)}{\partial d_s} + \dfrac{1}{2} \cdot \dfrac{\partial A_s(d_s, h_s)}{\partial h_s} < 0$, $\forall s \in D$.

**Remark:** According to the equivalence between (VI 4) and the combined equilibrium issue, it follows that if the assumptions in Theorem 4 are valid, then the combined equilibrium must be unique.

Let us consider some practical implications of the assumptions in Theorem 4. Firstly, suppose that $\frac{\partial A_s(d_s, h_s)}{\partial h_s} = \frac{\partial B_s(d_s, h_s)}{\partial d_s}$ is valid for any $s \in D$. That is to say, in business center $s$, consumers' marginal shopping utility with respect to $h_s$ equals to companies' marginal location utility with respect to $d_s$. This simplification can enable us to derive the practical implications of the assumptions in Theorem 4. Then the first assumption in Theorem 4 can be rewritten as

$$\frac{\partial A_s(d_s, h_s)}{\partial d_s} + \frac{\partial A_s(d_s, h_s)}{\partial h_s} < 0, \quad \forall s \in D. \tag{25}$$

From consumers' point of view, the inequality in (25) can be explained that the value of marginal congestion disutility $-\frac{\partial A_s(d_s, h_s)}{\partial d_s}$ is larger than the value of marginal shopping utility $\frac{\partial A_s(d_s, h_s)}{\partial h_s}$. In other words, the consumers think the additional congestion cost of increasing one more unit of consumers in $s$ exceeds the additional benefit of increasing one more unit of companies in $s$. On the other hand,

$$A_s(d_s+1, h_s+1) - A_s(d_s, h_s) \approx \frac{\partial A_s(d_s, h_s)}{\partial d_s}(d_s+1-d_s) + \frac{\partial A_s(d_s, h_s)}{\partial h_s}(h_s+1-h_s)$$

$$= \frac{\partial A_s(d_s, h_s)}{\partial d_s} + \frac{\partial A_s(d_s, h_s)}{\partial h_s}.$$

Thus the inequality in (25) can also be regarded as that the attraction of business center $s$ to consumers would decrease if one more unit of consumers and one more unit of business companies enter $s$ in the same time. In addition, the second assumption in Theorem 4 can be rewritten as

$$\frac{\partial B_s(d_s, h_s)}{\partial h_s} + \frac{\partial B_s(d_s, h_s)}{\partial d_s} < 0, \quad \forall s \in D. \tag{26}$$

Similarly, the inequality in (26) is to say that the business companies think the additional land rent cost of increasing one more unit of companies in $s$ exceeds the additional benefit of increasing one more unit of consumers in $s$, or the attraction of $s$ to companies would decrease if one more unit of consumers and one more unit of companies enter $s$ in the same time. From previous discussion, the prerequisites in Theorem 4 may demonstrate some cases in practice. They can be valid for the very congested business center in which there are numerous consumers and companies. In such business center, the marginal disutility of aggregation effect exceeds the marginal utility of aggregation effect for not only consumers but also companies.

## 3. The congestion pricing principles based on the combined equilibrium

The objective of this section is to investigate the congestion pricing principles based on the combined equilibrium. This section is divided into two parts. The first part aims at deriving the optimal road pricing scheme that can support a parametric traffic equilibrium as a "uniform" traffic system optimum. This means that the total social cost of travelers can achieve minimization. The second part is devoted to studying the general optimal congestion pricing scheme which consists of road pricing and business congestion tax. The goal of the scheme is to minimize the total social cost of travelers and business companies.

3.1 The optimal road pricing scheme

Now we are going to derive the optimal road pricing scheme that can minimize the total social cost of travelers. According to Yang's work (Yang, 1999), the total social cost of travelers in this paper can be measured as the total cost of travel time minus the travelers' benefit

$$\frac{1}{\alpha}\sum_{r\in O}\sum_{s\in D}\sum_{k\in R_{rs}} f_k^{rs}(\ln f_k^{rs} -1) + (\frac{1}{\gamma}-\frac{1}{\alpha})\sum_{r\in O}\sum_{s\in D} q_{rs}(\ln q_{rs} -1) + \sum_{a\in A} x_a t_a(x_a) - \sum_{s\in D} A_s(d_s, h_s)d_s ,$$

where $h = (h_s, s \in D)^T$ is the parametric vector.

Clearly, the parametric traffic system optimum requires the total social cost of travelers to be minimized based on current business location distribution $h = (h_s, s \in D)^T$. So we have the following parametric optimization problem.

**(OP 5)**

$$\min_{f,x,q,d} \quad \frac{1}{\alpha} \sum_{r \in O} \sum_{s \in D} \sum_{k \in R_{rs}} f_k^{rs} (\ln f_k^{rs} - 1) + (\frac{1}{\gamma} - \frac{1}{\alpha}) \sum_{r \in O} \sum_{s \in D} q_{rs} (\ln q_{rs} - 1) + \sum_{a \in A} x_a t_a(x_a) - \sum_{s \in D} A_s(d_s, h_s) d_s$$

s.t. $(f, x, q, d) \in K_1$

The above optimization model reflects the government's expectation that the travelers in the network consciously choose their destinations and travel routes to minimize their total social cost based on $h = (h_s, s \in D)^T$. The Kuhn-Tucker conditions for this optimization problem are

$$q_{rs} = O_r \frac{\exp\{-\gamma[v_{rs} - \bar{A}_s(d_s, h_s)]\}}{\sum_{j \in D} \exp\{-\gamma[v_{rj} - \bar{A}_j(d_j, h_j)]\}}, \quad \forall r \in O, s \in D \tag{27}$$

$$f_k^{rs} = q_{rs} \frac{\exp\{-\alpha[\sum_{a \in A} \delta_{ak}^{rs} \bar{t}_a(x_a)]\}}{\sum_{j \in R_{rs}} \exp\{-\alpha[\sum_{a \in A} \delta_{aj}^{rs} \bar{t}_a(x_a)]\}}, \quad \forall k \in R_{rs}, r \in O, s \in D \tag{28}$$

where $v_{rs}$ is the Lagrange multiplier and $v_{rs} = -\frac{1}{\alpha} \sum_{j \in R_{rs}} \exp\{-\alpha[\sum_{a \in A} \delta_{aj}^{rs} \bar{t}_a(x_a)]\}$,

$$\bar{t}_a(x_a) = t_a(x_a) + \hat{t}_a(x_a), \quad \bar{A}_s(d_s, h_s) = A_s(d_s, h_s) + \hat{A}_s(d_s, h_s), \tag{29}$$

where $\hat{A}_s(d_s, h_s) = d_s \frac{\partial A_s(d_s, h_s)}{\partial d_s}$, $\hat{t}_a(x_a) = x_a \frac{dt_a(x_a)}{dx_a}$.

According to $\frac{dt_a(x_a)}{dx_a} > 0$ and $\frac{\partial A_s(d_s, h_s)}{\partial d_s} < 0$, it follows that

$\hat{t}_a(x_a) \geq 0$, $\hat{A}_s(d_s, h_s) \leq 0$,

$$\bar{t}_a(x_a) = t_a(x_a) + \hat{t}_a(x_a) \geq t_a(x_a), \tag{30}$$

$$\bar{A}_s(d_s, h_s) = A_s(d_s, h_s) + \hat{A}_s(d_s, h_s) \leq A_s(d_s, h_s). \tag{31}$$

Obviously, $t_a(x_a)$ and $A_s(d_s, h_s)$ are respectively the actual travel time of the link $a$ and the actual attraction measure of the business center $s$ for travelers. From the point of view of travelers, the travel time is disutility and the destination attraction represents utility. Thus, from (30) and (31), $\hat{A}_s(d_s, h_s)$ and $\hat{t}_a(x_a)$ can be viewed as the negative externalities which stem from the effect of congestion on the links and the business centers in the network. That is to say, $\hat{A}_s(d_s, h_s)$ is the additional congestion disutility that a traveler imposes on all other travelers in $s$ while $\hat{t}_a(x_a)$ is the additional travel time that a traveler imposes on all other travelers in $a$. Furthermore, note that the Kuhn-Tucker conditions for (OP 5), namely (27) and (28), include the negative externalities and have the same forms as the parametric traffic equilibrium conditions if the link travel time function and the destination attraction function for travelers in (1) and (2) are replaced by $\bar{t}_a(x_a)$ and $\bar{A}_s(d_s, h_s)$ in (29). Therefore, by imposing a toll that exactly equals to the corresponding negative externality on each link and each business center's entrance, we can ensure that the travelers' optimal private travel choices, which include destination choices and route choices, will also be "relatively optimal" choices according to the minimization of total social cost of travelers. The word "relatively" refers to that the optimality is based on the current business location distribution $h = (h_s, s \in D)^T$. Evidently, the "relatively optimal" road pricing scheme is not qualified because this scheme is still dependent on the business location distribution.

Now the aim is to derive the uniform optimal scheme, which is independent of the business location distribution. Firstly, we should verify how travelers and companies would act respectively if the desired road pricing policy is carried out. Before establishing the mathematical model, we need to prove that (OP 5) is equivalent to its stationary point problem, which is given below.

**(VI 6)**

Find a vector $(f^*, x^*, q^*, d^*) \in K_1$ such that

$$\sum_{r \in O} \sum_{s \in D} \sum_{k \in R_{rs}} (f_k^{rs} - f_k^{rs*}) \frac{1}{\alpha} \ln f_k^{rs*} + \sum_{a \in A} (x_a - x_a^*)[t_a(x_a^*) + x_a^* \frac{dt_a(x_a^*)}{dx_a}]$$

$$+ \sum_{r \in O} \sum_{s \in D} (q_{rs} - q_{rs}^*)(\frac{1}{\gamma} - \frac{1}{\alpha}) \ln q_{rs}^* - \sum_{s \in D} (d_s - d_s^*)[A_s(d_s^*, h_s) + d_s^* \frac{\partial A_s(d_s^*, h_s)}{\partial d_s}] \geq 0,$$

$$\forall (f, x, q, d) \in K_1.$$

Note that there is a parametric vector $h = (h_s, s \in D)^T$ in (VI 6), which indicates the influence of the location distribution of business companies. When only the road pricing policy is allowed to carry out, we can not expect the companies to consciously act according to the will of the government since the road pricing policy has no direct influence on them. The business companies' location distribution is still based on the principle of personal utility maximization. Thus the behaviors of companies can still be described by (VI 2) under the policy of road pricing. By the equivalence between (OP 5) and (VI 6), it follows that now the behaviors of travelers and companies can be precisely depicted by the following system of variational inequalities, which is a combination of (VI 2) and (VI 6).

**(VIS 7)**

Find $(f^*, x^*, q^*, d^*) \in K_1$ and $h^* \in K_2$ such that

$$\sum_{r \in O} \sum_{s \in D} \sum_{k \in R_{rs}} (f_k^{rs} - f_k^{rs*}) \cdot \frac{1}{\alpha} \ln f_k^{rs*} + \sum_{a \in A} (x_a - x_a^*)[t_a(x_a^*) + x_a^* \frac{dt_a(x_a^*)}{dx_a}]$$

$$+ \sum_{r \in O} \sum_{s \in D} (q_{rs} - q_{rs}^*) \cdot (\frac{1}{\gamma} - \frac{1}{\alpha}) \ln q_{rs}^* - \sum_{s \in D} (d_s - d_s^*)[A_s(d_s^*, h_s^*) + d_s^* \frac{\partial A_s(d_s^*, h_s^*)}{\partial d_s}] \geq 0,$$

$$\sum_{s \in D} (h_s - h_s^*)[\frac{1}{\beta} \ln h_s^* - B_s(d_s^*, h_s^*)] \geq 0, \qquad \forall (f, x, q, d) \in K_1, h \in K_2.$$

It is easy to prove that the Kuhn-Tucker conditions of (VIS 7) are just (3), (27) and (28). Since the Kuhn-Tucker conditions for a variational inequality with linear constraints is equivalent to the variational inequality itself (Facchinei and Pang, 2003), the equations (27), (28) and (3) give a complete and accurate characterization of the behaviors of travelers and companies under the road pricing policy. Therefore, by comparing (27) with (1) and (28) with (2) respectively, we can obtain the following uniform optimal road pricing scheme which can ensure that the total social cost of travelers achieves minimization over all possible business location distribution patterns in the state of equilibrium.

**The uniform optimal road pricing scheme**

(a) Find the global optimal solution $(f, x, q, d)$ of following optimization problem:

$$\min_{f,x,q,d} \frac{1}{\alpha} \sum_{r \in O} \sum_{s \in D} \sum_{k \in R_{rs}} f_k^{rs} (\ln f_k^{rs} - 1) + (\frac{1}{\gamma} - \frac{1}{\alpha}) \sum_{r \in O} \sum_{s \in D} q_{rs} (\ln q_{rs} - 1) + \sum_{a \in A} x_a t_a(x_a) - \sum_{s \in D} A_s(d_s, h_s) d_s$$

s.t.  $(f, x, q, d, h)$ belongs to the solution set of (VIS 7).

(b) For any link $a \in A$, any traveler using this link will be charged a fee that equals to $x_a \frac{dt_a(x_a)}{dx_a}$.

(c) For any $s \in D$, any traveler entering this business center will be charged a fee that equals to $-d_s \frac{\partial A_s(d_s, h_s)}{\partial d_s}$.

**Remark 1:** Step (a) is to select the best equilibrium solution that can uniformly minimize the travelers' total social cost from all possible equilibrium solutions. If the solution of (VIS 7) is unique, then the step (a) would be is simplified as "Find the unique solution of (VIS 7)".

**Remark 2:** Since $A_s(d_s, h_s)$ is strictly decreasing with $d_s$, $-d_s \frac{\partial A_s(d_s, h_s)}{\partial d_s} \geq 0$. So the step (c) is to charge the travelers, not to provide them with subsidies.

The traditional link-based optimal road pricing scheme is just the subscheme (b) of the uniform optimal road pricing scheme. The differences between the traditional scheme and the uniform optimal scheme result from the fact that the traditional scheme assumes the business location distribution pattern is given in advance without considering the important influence of the business location distribution. Hence, within the theoretic framework in this paper, the traditional scheme can not derive a combined equilibrium toward a uniform traffic system optimum or even a "relative" traffic system optimum which depends on the business location distribution (the "relative" traffic system optimum is depicted by the parametric optimization problem (OP 5)). Thus the uniform optimal road pricing scheme generalizes the traditional link-based optimal road pricing scheme.

3.2 The optimal congestion pricing scheme that consists of road pricing and business congestion tax

The objective of this part is to derive the optimal congestion pricing scheme that contains road pricing and business congestion tax. The scheme should consider the social cost of both travelers and companies because the business congestion tax has a direct influence on the income of companies and the business location distribution pattern. At first, we would like to point out the rationality of the business congestion tax policy which aims at the group of business companies. It can be asserted that the business companies should also be partly responsible for the congestion since the aggregation of business companies may influence the trip distribution of travelers. In fact, from a mathematical point of view, by (1) and (2), we have:

$$f_k^{rs} = O_r \frac{\exp\{-\alpha[v_{rs} - A_s(d_s, h_s)]\}}{\sum_{j \in D} \exp\{-\alpha[v_{rj} - A_j(d_j, h_j)]\}} \cdot \frac{\exp\{-\alpha[\sum_{a \in A} \delta_{ak}^{rs} t_a(x_a)]\}}{\sum_{j \in R_{rs}} \exp\{-\alpha[\sum_{a \in A} \delta_{aj}^{rs} t_a(x_a)]\}},$$

$$\forall k \in R_{rs}, r \in O, s \in D.$$

Clearly, the above equality indicates that the traffic flow on every path in the network depends on the location distribution of business companies. Therefore, the traditional road pricing scheme may be unfair to the travelers because the scheme only charges them. Thus the business congestion tax policy is necessary and reasonable.

Now we start to establish the model for the overall system optimum under the policy of road pricing and business congestion tax. The model should reflect the principle on which the consumers and business companies jointly insist, as we expect, to minimize the total social cost of travelers and companies. Based on this idea, we present the following optimization model.

**(OP 8)**

$$\min_{f,x,q,d,h} \frac{1}{\alpha}\sum_{r\in O}\sum_{s\in D}\sum_{k\in R_{rs}} f_k^{rs}(\ln f_k^{rs}-1) + (\frac{1}{\gamma}-\frac{1}{\alpha})\sum_{r\in O}\sum_{s\in D} q_{rs}(\ln q_{rs}-1) + \frac{1}{\beta}\sum_{s\in D} h_s(\ln h_s-1) + \sum_{a\in A} x_a t_a(x_a)$$

$$-\sum_{s\in D} A_s(d_s,h_s)d_s - \sum_{s\in D} B_s(d_s,h_s)h_s$$

s.t. $(f,x,q,d,h) \in K_1 \times K_2$

The economic implication of (OP 8) is similar to that of (OP 5). Now we want to show that solving the optimization problem (OP 8) leads to the desired optimal congestion pricing scheme. It is easy to prove that the Kuhn-Tucker conditions of (OP 8) are

$$q_{rs} = O_r \frac{\exp\{-\gamma[v_{rs}-\bar{A}_s(d_s,h_s)]\}}{\sum_{j\in D}\exp\{-\gamma[v_{rj}-\bar{A}_j(d_j,h_j)]\}}, \quad \forall r\in O, s\in D, \tag{32}$$

$$f_k^{rs} = q_{rs} \frac{\exp\{-\alpha[\sum_{a\in A}\delta_{ak}^{rs}\bar{t}_a(x_a)]\}}{\sum_{j\in R_{rs}}\exp\{-\alpha[\sum_{a\in A}\delta_{aj}^{rs}\bar{t}_a(x_a)]\}}, \quad \forall k\in R_{rs}, r\in O, s\in D, \tag{33}$$

$$h_s = T \frac{\exp[\beta\bar{B}_s(d_s,h_s)]}{\sum_{j\in D}\exp[\beta\bar{B}_j(d_j,h_j)]}, \quad \forall s\in D, \tag{34}$$

where $v_{rs}$ is the Lagrange multiplier and $v_{rs} = -\frac{1}{\alpha}\ln\{\sum_{j\in R_{rs}}\exp[-\alpha(\sum_{a\in A}\delta_{aj}^{rs}\bar{t}_a(x_a))]\}$,

$$\bar{A}_s(d_s,h_s) = A_s(d_s,h_s) + d_s\frac{\partial A_s(d_s,h_s)}{\partial d_s} + h_s\frac{\partial B_s(d_s,h_s)}{\partial d_s}, \tag{35}$$

$$\bar{t}_a(x_a) = t_a(x_a) + x_a\frac{dt_a(x_a)}{dx_a}, \tag{36}$$

$$\bar{B}_s(d_s, h_s) = B_s(d_s, h_s) + d_s \frac{\partial A_s(d_s, h_s)}{\partial h_s} + h_s \frac{\partial B_s(d_s, h_s)}{\partial h_s}. \tag{37}$$

In above equations, $A_s(d_s, h_s)$ and $B_s(d_s, h_s)$ are respectively the actual attraction measure for travelers and companies in the business center $s$. The economic implications of the other terms in (35), (36) and (37) are given below.

1. $d_s \dfrac{\partial A_s(d_s, h_s)}{\partial d_s}$ is the additional congestion disutility that a traveler imposes on all other travelers in the business center $s$. (38)

2. $h_s \dfrac{\partial B_s(d_s, h_s)}{\partial d_s}$ is the additional business benefit that a traveler brings to the companies in the business center $s$. (39)

3. $x_a \dfrac{dt_a(x_a)}{dx_a}$ is the additional travel time that a traveler imposes on all other travelers in the link $a$. (40)

4. $d_s \dfrac{\partial A_s(d_s, h_s)}{\partial h_s}$ is the additional consuming utility that a company brings to the travelers who choose the business center $s$ as their destination. (41)

5. $h_s \dfrac{\partial B_s(d_s, h_s)}{\partial h_s}$ is the additional land rent cost that a company imposes on all other companies in the business center $s$. (42)

Therefore, the terms in (39)-(42) can be treated as the externalities due to the aggregation effect of travelers and companies in the network. To internalize the externalities, a toll should be imposed on each link and each business center's entrance and it should equals to the corresponding externality. To be precise, we obtain the following pricing scheme.

**The optimal congestion pricing scheme**

(d) Find the global optimal solution $(f, x, q, d, h)$ of the optimization problem (OP 8).

(e) For any link $a \in A$, any traveler using this link is charged a fee that equals to $x_a \dfrac{dt_a(x_a)}{dx_a}$.

(f) For any $s \in D$, any traveler entering this business center is charged a fee that equals to $u_s = -[d_s \dfrac{\partial A_s(d_s, h_s)}{\partial d_s} + h_s \dfrac{\partial B_s(d_s, h_s)}{\partial d_s}]$.

(g) For any $s \in D$, any company in this business center is imposed a tax that amounts to $v_s = -[d_s \dfrac{\partial A_s(d_s, h_s)}{\partial h_s} + h_s \dfrac{\partial B_s(d_s, h_s)}{\partial h_s}]$.

In general, it is difficult to identify that $u_s$ and $v_s$ are positive or negative. Therefore, $u_s > 0$ means that the travelers are charged a fee, $u_s < 0$ means that the travelers are provided with a subsidy, $v_s > 0$ indicates that the companies are imposed an additional tax, and $v_s < 0$ indicates that the companies are given an allowance or a tax break. The possible rebate policy in Step (f) can be implemented in such a way that every traveler who has shopped in the specified business centers can receive the subsidy through providing his/her shopping list to the government.

Now we further study the economic meaning of the optimal congestion pricing scheme. We are especially interested in that of the steps (f) and (g). Firstly, according to (38)-(42), $d_s \dfrac{\partial A_s(d_s, h_s)}{\partial d_s}$ and $d_s \dfrac{\partial A_s(d_s, h_s)}{\partial h_s}$ are respectively the negative externality and the positive externality to travelers while $h_s \dfrac{\partial B_s(d_s, h_s)}{\partial d_s}$ and $h_s \dfrac{\partial B_s(d_s, h_s)}{\partial h_s}$ are respectively the positive externality and the negative externality to companies. Secondly, the externalities $d_s \dfrac{\partial A_s(d_s, h_s)}{\partial d_s}$ and $h_s \dfrac{\partial B_s(d_s, h_s)}{\partial d_s}$ are caused by the aggregation of travelers, and the externalities $d_s \frac{\partial A_s(d_s,h_s)}{\partial h_s}$ and $h_s \frac{\partial B_s(d_s,h_s)}{\partial h_s}$ are due to the aggregation of companies. Therefore, in the business center $s$, every individual (a traveler or a company) may simultaneously impose a positive external effect and a negative external effect on the overall system.

$$-u_s = d_s \frac{\partial A_s(d_s,h_s)}{\partial d_s} + h_s \frac{\partial B_s(d_s,h_s)}{\partial d_s} \quad \text{and} \quad -v_s = d_s \frac{\partial A_s(d_s,h_s)}{\partial h_s} + h_s \frac{\partial B_s(d_s,h_s)}{\partial h_s}$$

are respectively the overall external effect of a traveler and a company which are imposed on the system. If $-u_s > 0$, it can be asserted that every traveler in $s$ has a positive overall external effect and a rebate policy for the travelers in this center should be introduced to compensate the positive externalities. On the other hand, if $-u_s < 0$, it indicates that every traveler in $s$ has a negative overall external effect and a charge policy is needed to eliminate the negative externalities. A similar discussion about the externalities caused by the aggregation of companies in every business center can be taken according to the sign of $-v_s$. All in all, the core of the optimal congestion pricing scheme which includes road pricing and business congestion tax can be summarized as the following sentence.

*For every individual, if his private choice has a positive overall external effect on the system, he will receive a subsidy, otherwise, he will have to pay a toll, and the subsidy/toll should be equivalent to the external effect.*

**4. A numerical example**

In this section, we present a numerical example to demonstrate that the optimal congestion pricing scheme, which is proposed in Section 3.2, can indeed reduce the social cost. We primarily concern how many social costs of overall system can be reduced under the pricing scheme, and how the scheme influences the travelers' trip distribution and the companies' location distribution. The road network, as shown in Fig 4.1, has 6 nodes and 7 links, where node 1 and node 2 are origin nodes, and node 4 and node 5 are destination nodes (business centers).

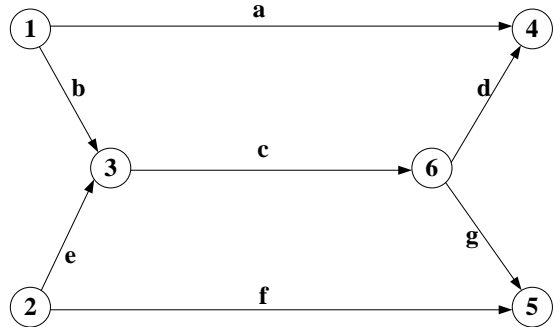

**Fig. 4.1.**

The numerical experiments are implemented for twelve scenarios indexed by $k = 1,2,...,12$, varying the travel demands from every origin node according to $O_1^1 = 40, O_2^1 = 60$, $O_s^k = O_s^{k-1} + 5$ $(k = 2,3,...,12; s = 1,2)$, and keeping other factors unchanged. The total number of companies in the network, $T$, is set to be 50. The numerical results are plotted in the following figures.

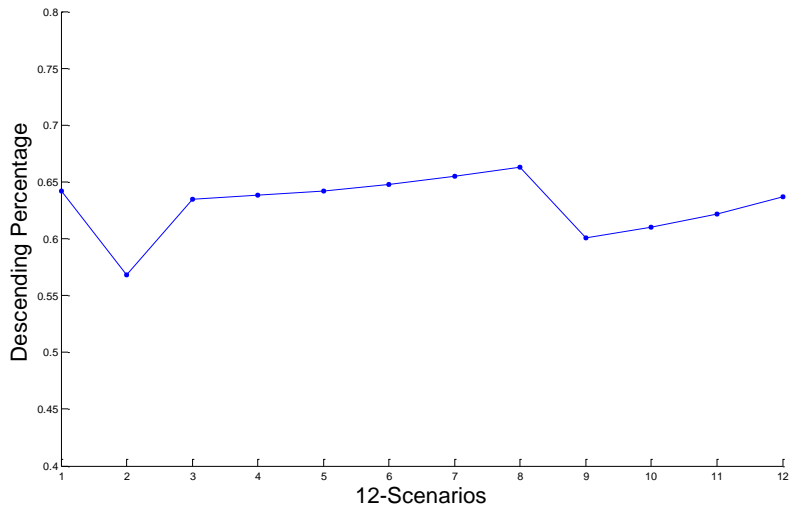

Fig. 4.2. The descending percentage of total social cost.

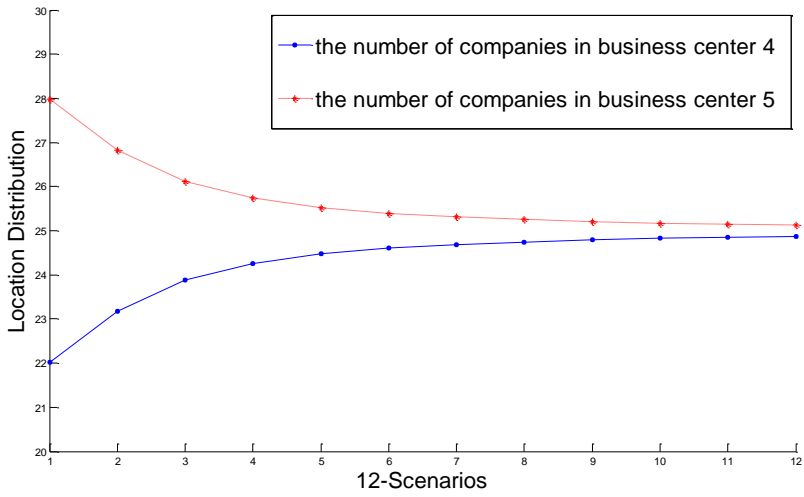
Fig. 4.3. Companies' location distribution with toll charge.

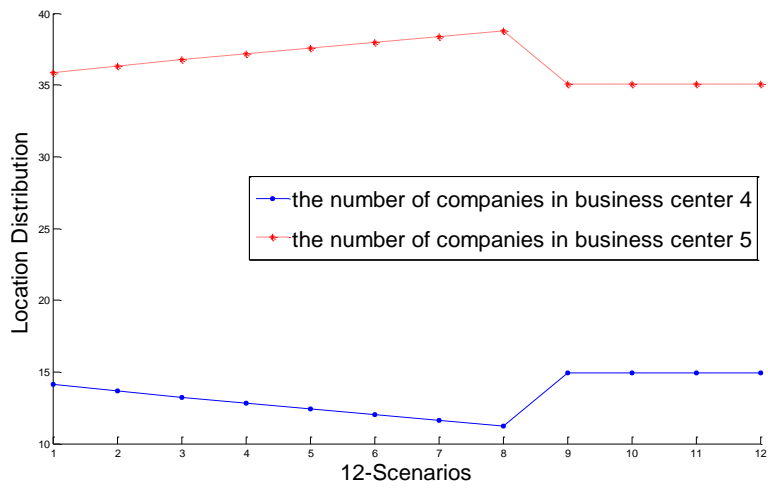
Fig. 4.4. Companies' location distribution with no toll charge.

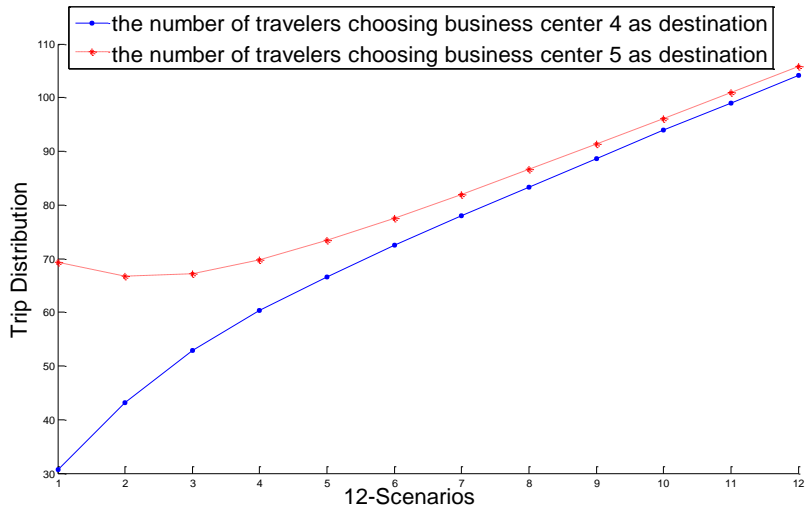
Fig. 4.5. Travelers' trip distribution with toll charge.

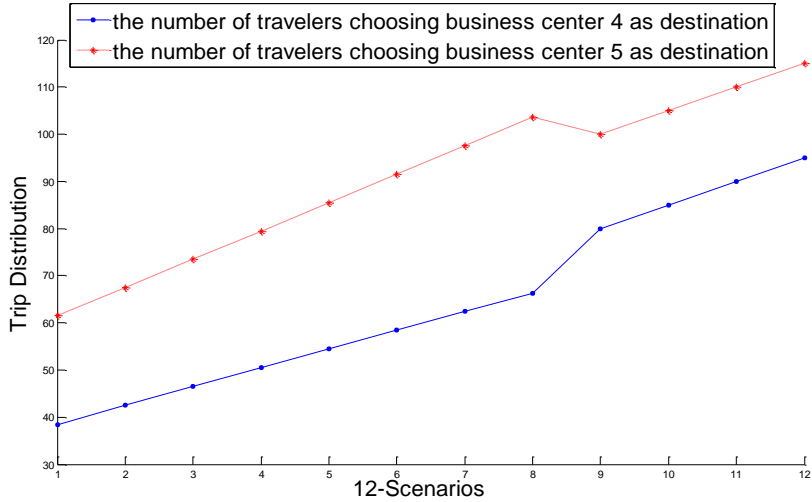
Fig. 4.6. Travelers' trip distribution with no toll charge.

It can be seen from Fig. 4.2 that the social cost of overall system is reduced for every scenario under the optimal congestion pricing scheme. To be specific, the descending percentage of social cost of overall system remains close to 60%, which indicates that the pricing scheme may be efficient and stable in reducing total social cost. Furthermore, by comparing Fig. 4.3 with Fig. 4.4 and Fig. 4.5 with Fig. 4.6 respectively, it follows that travelers' trip distribution pattern and companies' location distribution pattern are more balanced under the optimal congestion pricing policy.

## 5. Conclusion

In this working paper, we propose a new combined equilibrium model of business land-use and investigate its congestion pricing principles. Firstly, the model integrates travelers' traffic equilibrium with business companies' competitive location equilibrium, which can be depicted by two parametric variational inequalities respectively. Then by combing them together, we obtain a variational inequality equivalent to the combined equilibrium. The existence and uniqueness of equilibrium solution are established. Additionally, the mathematical and economic principles of congestion pricing associated with the combined equilibrium are studied. Mathematically, we prove that there is an optimal road pricing scheme that can minimize the social cost of travelers. This scheme generalizes the traditional link-based optimal road pricing scheme and reduces the social cost of travelers within a more general framework. Furthermore, if allowed to simultaneously impose charges on travelers and companies, we prove that there exists an optimal congestion pricing scheme that can derive a combined equilibrium toward an overall system optimum according to the minimization of the total social cost of travelers and companies. The economic meaning of every pricing scheme proposed in this paper is investigated in detail. A numerical example is presented to demonstrate that the previous optimal congestion pricing scheme can indeed reduced the social cost.